\documentclass[12pt, twoside,en]{article}
\usepackage[utf8]{inputenc}

\usepackage{amsthm}
\usepackage{blkarray}
\usepackage{tikz-cd}
\usepackage{amsfonts} 
\usepackage{mathtools}
\usepackage{hyperref}
\usepackage{amssymb}
\usepackage{tikz}
\usepackage[affil-it]{authblk}
\usepackage{blkarray}
\usepackage[backend=biber]{biblatex}
\usepackage[toc,page]{appendix}
\usepackage[english]{babel}
\usepackage[ruled,vlined]{algorithm2e}
\usepackage{enumerate}

\AtBeginDocument{%
  \providecommand\BibTeX{{%
    \normalfont B\kern-0.5em{\scshape i\kern-0.25em b}\kern-0.8em\TeX}}}

\newtheorem{theorem}{Theorem}[section]
\newtheorem{proposition}{Proposition}[section]
\newtheorem{definition}{Definition}[section]
\newtheorem{lemma}{Lemma}[section]

\newtheorem{corollary}{Corollary}[section]

\newtheorem{que}{Question}[section]

\newcommand\restr[2]{{
  \left.\kern-\nulldelimiterspace 
  #1 
  \vphantom{\big|} 
  \right|_{#2} 
  }}
 
\usepackage[margin=1in]{geometry}

\DeclareMathSymbol{\ast}{\mathbin}{symbols}{"03}

\title{A Lower Bound for the Rank of Matroid Intersection
}

\author{Liu Tianyu}
\affil{Department of 
Mathematics\\ University of Chinese Academy of Sciences \\ China-100049}
\date{\vspace{-5ex}}

\begin{document}

\maketitle

\begin{abstract}
Matroid is a generalization of many fundamental objects in combinatorial mathematics , and matroid intersection problem is a classical subject in combinatorial optimization . However , only the intersection of two matroids are well understood  . The solution of the intersection problem of more than three matroids is proved to be \textbf{NP-hard} . We will give a lower bound  estimate on the maximal cardinality of the common independent sets in matroid intersections .   We will also study some properties of the intersection of more than two matroids and deduce some analogous results for Edmonds' Min-max theorems for matroids intersection .  
\end{abstract}

%


\section{Introduction}
    Matroid was firstly introduced by Hassler Whitney in 1935 \cite{whitney1935abstract} and also independently discovered by Takeo Nakasawa, whose work was forgotten for many years \cite{nishimura2009lost} . In 1950s and 1960s , T. Tutte made the prominent contributions to matroid theory in his outstanding papers ( Interested readers can check   \cite{tutte1959matroids} \cite{tutte1964lectures} )  . For a detailed introduction to matroid theory  , the readers are referred to \cite{neel2009matroids} \cite{oxley2006matroid} \cite{welsh2010matroid} . Other major contributors include Jack Edmonds, Jim Geelen, Eugene Lawler, László Lovász, Gian-Carlo Rota, P. D. Seymour, and Dominic Welsh . Without extra announcement a matroid is always meant to be a finite matroid . In this paper we will not
   discuss infinite matroid .

Many problems in combinatorial optimization can be reformulated as : given a system $(E,\mathcal{F})$ where \textit{E} is a finite set and $\mathcal{F}\subseteq 2^E$ , and a valued function $c:\mathcal{F}\to \mathbb{R}$ , we intend to find an element of $\mathcal{F}$ maximizing or minimizing $c$ . And many problems in combinatorial optimization can be restated as matroid problems . In this paper we restrict our discussion on matroids  .


Firstly let's recall Edmonds' min-max theorem for the intersection of two matroids , which plays an important role in Edmonds' algorithm
for matroid intersection .

\begin{theorem}[Edmonds,1970\cite{edm}]
\label{Edmonds}
  Let $(E,\mathcal{F}_1)$ and $(E,\mathcal{F}_2)$ two matroids , with rank functions $r_1$ and $r_2$ respectively . Then we have
  \[\max\{|X|:X\in\mathcal{F}_1 \cap\mathcal{F}_2\}=\min\{r_1(Q)+r_2(E\backslash Q):Q\subseteq E\}.\]
\end{theorem}
For intersections of more than two matroids , no analogous result of $\mathbf{Theorem\;1.}$ has been known to the author . And it has been proved that finding a subset with the maximal cardinality of the intersection of more than two matroids , is an $\mathbf{NP-hard}$ problem ( c.f.\cite{kor} Chapter 15 , Exercice 14 (c) ) . However ,  if we make some futher assumptions on given matroids (more than two) , we can formally generalize \ref{Edmonds} to the following .

\begin{proposition}
\label{prop1.1}
Let $(E,\mathcal{F}_1)$ , $(E,\mathcal{F}_2)$ , ... , $(E,\mathcal{F}_m)$ be m  matroids with rank functions $r_1,\dots,r_m$ respectively . And assume that $(E,\mathcal{F}_1\bigcap\mathcal{F}_2)$ , $(E,\mathcal{F}_1\bigcap\mathcal{F}_2\bigcap\mathcal{F}_3)$ , $\dots$ , $(E,\bigcap_{i=1}^m\mathcal{F}_i)$ are all matroids . And for any $X\subseteq E$  we take a increasing filtration of $X$ randomly :
\[X\subseteq X_1\subseteq X_2\subseteq \dots\subseteq X_{m-2}\subseteq E\]
Then we have the following min-max result :
\begin{equation}
    \begin{split}
        \max\{|S|:S\in\bigcap_{i=1}^m\mathcal{F}_i\}&=\min\{r_1(X)+r_2(X_1\backslash X)+\dots+r_{m-1}(X_{m-2}\backslash X_{m-3})+r_m(E\backslash X_{m-2}):\\ &\quad X\subseteq X_1\subseteq X_2\subseteq \dots\subseteq X_{m-2}\subseteq E\}
    \end{split}
\end{equation}
\end{proposition}

In fact the above theorem can be intuitively  obtained by inductions . We will give a proof in Section 3 . However ,
it seems that the assumption we put on the above proposition is too strong : the intersection $\bigcap_{i=1}^k \mathcal{M}_i$ are all matroids for $k=1,\dots,m$ . 

We have the following rough upper bound estimation for matroid intersection . 

\begin{corollary}
\label{cor1.1}
Let $\{\mathcal{M}_i(E,\mathcal{F}_i)\}_{i=1}^m$ be a set of $m$ matroids on $E$ and $r_i$ is the rank function of $\mathcal{M}_i$ , then we have 

\begin{equation}
    \max_{I\in\mathcal{B}(\bigcap_{i=1}^m\mathcal{M}_i)}|I|\leqslant \min \sum_{i=1}^m r_i(X_i)
\end{equation}
where the minimum in the right-hand-side  is taken over all the partition  $\bigsqcup_{i=1}^m X_i=E$ of $E$ .
\end{corollary}

The following result is fumdamental for studying the union of matroids over a common ground set and is useful in our proof of main results in this paper .

\begin{theorem}[Nash-Williams,1967\cite{whi}]
\label{thm1.2}
Let $\mathcal{M}_1=(E,\mathcal{F}_1)$ , $\mathcal{M}_2(E,\mathcal{F}_2)$ , ... , $\mathcal{M}_m(E,\mathcal{F}_m)$ be  matroids with rank functions $r_1,\dots,r_n$ respectively . And denote their (disjoint) union as $(E,\mathcal{F})=$ $(E,\bigsqcup_{i-1}^n \mathcal{F}_i)$ . Then $(E,\mathcal{F})$ is a matroid with rank function
\[r(X)=\min_{A\subseteq X} \biggl( |X\backslash A|+\sum_{i-1}^k r_i(A) \biggr) .\]
  
\end{theorem}

Using some properties of the duality of a matroid ,  we've proved the following main result of this paper which can give an lower bound  for the maximal cardinality of the intersection of n matroids :
\begin{proposition}
\label{prop1.2}
  Let $\mathcal{M}_1=(E,\mathcal{F}_1)$ , $\mathcal{M}_2=(E,\mathcal{F}_2)$ , $\dots$ , $\mathcal{M}_m=(E,\mathcal{F}_m)$ be m matroids on the ground set E , with rank functions $r_1 ,\dots$ $r_m$ respectively , here E is a finite set . And set $(E,\mathcal{F}_j^{*})$ the duality of $(E,\mathcal{F}_j)$ for $j=1,\dots,m$ . Moreover we assume that $(E,\bigcap_{i=1}^m\mathcal{F}_i)$ is an independent system with rank function $\Bar{r}$ ( Notice that the intersection of matroids are not generally a matroid ) . Then we have 
  \[\Bar{r}^*(X)\leqslant \min_{A \subseteq E}\biggl(|X|+(n-1)|A|+\sum_{i=1}^m r_i(E\backslash A)-\sum_{i=1}^m r_i(E)\biggr)\]
    Here $\Bar{r}^*$ is the rank function of $(E,(\bigcap_{i=1}^m\mathcal{F}_i)^*)$ . Especially if we take $X=E$ , we have the following :
  \[\max\{|X|:X\in \bigcap_{i=1}^m\mathcal{F}_i\}\geqslant \max_{A \subseteq E}\biggl(\sum_{i=1}^m r_i(E)-\sum_{i=1}^m r_i(E\backslash A)-(m-1)|A|\biggr).\]
  
If we restrict ourselves on two matroids  $(E,\mathcal{F}_1),(E,\mathcal{F}_2)$ , and furthermore assume that $\mathcal{F}_1\bigcap \mathcal{F}_2$ is a matroid. The two  inequalities above are actually equalities .
\end{proposition}

In fact, the lower bound of the maximum cardinality in the matroid intersection  we deduce above is too weak to be interesting. We will give a stronger lower bound estimation for the maximum cardinality in the following theorem. In this paper, we use $\mathcal{B}(\mathcal{M})$ to denote the set of bases in a given independent system $\mathcal{M}$
 .

\begin{theorem}
\label{thm1.3}

We set $\{T_i\}_{i=1}^m$ is a covering and co-covering of $E$ , here covering means that $\bigcup_{i=1}^m T_i =E$ and co-covering means $\bigcup_{i=1}^m T_i^c =E$ , where $T_i^c = E-T_i$ is the complement of $T_i$ for $i=1.\dots,m$ . And let $\mathcal{M}_1=(E,\mathcal{F}_1)$ , $\mathcal{M}_2=(E,\mathcal{F}_2)$ , $\dots$ , $\mathcal{M}_m=(E,\mathcal{F}_m)$ be m matroids as above . Moreover we assume that $T_i\in \mathcal{B}(\mathcal{M}_i)$ resp. for all $i$  , then we have :
\begin{equation}
    \max_{S\in\mathcal{B}( \bigcap_{i=1}^m\mathcal{F}_i)}|S|\geqslant \max_{A\subseteq E}\biggl(2\cdot\sum_{i=1}^mr_i(E)-(m-1)|A|-\sum_{i=1}^m r_i(E-A \cap T_i^c)\biggr)
\end{equation}

\end{theorem}

 \textbf{Organization.}  The remaining part of this paper is organized as follow : In section 2 we will recall the definition of matroid and list their basic properties without proofs . In section 3 we show the relations of the maximal cardinality of the independent sets of a matroid and its duality , using these dual relations we can give a proof of \textbf{Proposition}\ref{prop1.2} and \textbf{Theorem}\ref{thm1.3} .

\section{Some Preliminaries for Matroid Theory}

\subsection{Basic Properties of Matroids}

 In this subsection we will give the definition of matroids and list some basic facts without proofs just give the standard references for interested readers .

 Almost all notations in this paper are followed as \cite{kor} except that we will use a different notation for the union of matroids and introduce some new concepts . In this paper we use $"\bigsqcup"$ to denote disjoint union of sets and union of matroids .
 
 \begin{definition}
 A system of a set is a pair $(E,\mathcal{F})$ where $E$ is a finite set and $\mathcal{F}\subseteq 2^E$ . And an $\mathbf{matroid}$ is a system of set that satisfying the following three axioms:\\
 (M1) $\emptyset\in\mathcal{F}$;\\
 (M2)If $X\subseteq Y\in \mathcal{F}$ , $X\in \mathcal{F}$ ;\\
 (M3)If $X,Y\in \mathcal{F}$ and $|X|>|Y|$ , then there exists an element $x\in X\setminus Y$ such that $Y\bigcup \{x\}\in \mathcal{F}$ .

 \end{definition}
 
 \begin{definition}
 The $\mathbf{rank}$ of a matroid is $r(X):=\max \{|Y|:Y\subseteq X,Y\in \mathcal{F}\}$ , for any $X\subseteq E$ . 
 \end{definition}
 
 The following properties of rank function of a matroid can be straightforwardly obtained by definitions :
 
 \begin{proposition}
  Set $(E,\mathcal{F})$ a matroid with rank function r , then we have : for all $X,Y\subseteq E$\\
 $(R1)$ $r(X)\leq |X|$ ;\\
 $(R2)$ If $X\subset Y$ , then $r(X)\leq r(Y)$ ;\\
 $(R3)[\mathbf{submodularity}]$ $r(X\bigcup Y)+r(X\bigcap Y)\leq r(X)+r(Y)$ . 
 \end{proposition}

 To make things down to earth we give some examples of matroids in the following propositions . And it's easy to check that the following three systems satisfy the three axioms of matroids . 
 
 \begin{proposition}
 The following systems are all matroids:\\
 $(a)$ $A$ is a matrix over a field $\mathbb{F}$ , $E$ the set consisting of columns of $A$ , and $\mathcal{F}$:=$\{$ $F\subseteq E$ : the columns in F is linearly independent over $\mathbb{F}$ $\}$ .\\
 $(b)$ $E$ is the set consisting all edges of a given undirected graph $G$ , $\mathcal{F}:=\{F\subseteq E:(V(G),F) \}$ is a forest .\\
 $(c)$ Let $E$ be a finite set , $k$ be an integer , $\mathcal{F}:=\{F\subseteq E:|F|\leq k\}$ .
 \end{proposition}
    
    Another basic notion in matroid theory is the duality of a matroid which is important in our proof of main result in this paper .
    
    \begin{definition}
    Let $(E,\mathcal{F})$ be a independent system (i.e. , systems satisfying merely axioms M1 and M2 in the definition of matroids .) . Then the duality of $(E,\mathcal{F})$ is denoted as $(E,\mathcal{F*})$ , here:\\
    \[\mathcal{F^*}=\{F\subseteq E :There\; is\;a\;  basis\; B\; of\; (E,\mathcal{F})\; such\; that\; F\cap B=\emptyset\}.\]
    
    \end{definition}
    
It's clear that the duality of an independent system is also an independent system .

\begin{proposition}
\label{dualdual}
$(E,\mathcal{F^{**}})=(E,\mathcal{F}).$
\end{proposition}

\begin{proof}
$F\in \mathcal{F^{**}}\iff$ there is a basis $B^*$ of $(E,\mathcal{F^*})$ such that $F\bigcap B^*=\emptyset \iff$ there is a basis $B$ of $(E,\mathcal{F})$ s.t. $F\bigcap (E\setminus B)=\emptyset\iff$ $F\in \mathcal{F}$ .
\end{proof}

We will introduce two basic but important constrcutions called the \textbf{deletion} and \textbf{contraction} of a given matroid $\mathcal{M}$ .

\begin{definition}[Deletion]
Given a matroid $\mathcal{M}=(E,\mathcal{F})$ and a subset X of E , the deletion $\mathcal{M\setminus X}$ is a matroid $(E\setminus X , \mathcal{F}(\mathcal{M}\setminus X))$ , where 
\[\mathcal{F}(\mathcal{M}\setminus X)=\{C\subseteq E-X:C\in\mathcal{F}\}\]
.
\end{definition}
 
 It is easy to check that $\mathcal{F}(\mathcal{M}\setminus X)$ satisfies three matroid axioms .
 
 \begin{definition}[Contraction]
 Let $\mathcal{M}$ and X be as above , then the contraction of X in $\mathcal{M}$ is defined as 
 \[\mathcal{M}/X:=(\mathcal{M}^*\setminus X)^*\]
 .
 \end{definition}
 
 \begin{proposition}
 
 \begin{enumerate}[i)]
     \item $(\mathcal{M}\setminus X)^*=\mathcal{M}^*/X$,
     \item $(\mathcal{M} /X)^*=\mathcal{M}^*\setminus X$,
     \item \label{prop2.4.3} $r_{\mathcal{M}/X}(A)=r_{\mathcal{M}}(A\bigcup X)-r_{\mathcal{M}}(X)$ .
     
 \end{enumerate}
 
 \end{proposition}

  The following proposition ( especially property $(ii)$ ) is useful in our proof .
  
\begin{proposition}
\label{dualrank}
Let $(E,\mathcal{F})$ be an independent set , $(E,\mathcal{F}^*)$ be its duality . And put $r$ and $r^*$ be the rank function of $(E,\mathcal{F})$ and $(E,\mathcal{F}^*)$ respectively . Then we have \\
$(i)$ $(E,\mathcal{F})$ is a matroid $i.f.f.$ $(E,\mathcal{F}^*)$ is a matroid .\\
$(ii)$ If $(E,\mathcal{F})$ is a matroid , then the following equality
\[r^*(\mathcal{F})=|F|+r(E\setminus F)-r(E)\]
holds for all $F\subseteq E$ .
\end{proposition}

Finally in this section we will define the union of matroids over a common ground set $E$ . Set $(E,\mathcal{F}_1),\dots,(E,\mathcal{F}_n)$ be n matroids and
$(E,\mathcal{F}_1^*),\dots,(E,\mathcal{F}_n^*)$ there dualities respectively . We call a set $X\subseteq E$ is $partible$ if there is a partition $X=X_1\bigsqcup\dots\bigsqcup X_n$ , such that $X_i\in \mathcal{F}_i$ holds for every $i=1,\dots,n$ . And let $\bigsqcup_{i=1}^n \mathcal{F}_i$ be the family consisting of all the partible sets in $E$ associated to $(E,\mathcal{F}_1),\dots,(E,\mathcal{F}_n)$ .
Then we call $\bigsqcup_{i=1}^n \mathcal{F}_i$ the $\mathbf{union}$  or $\mathbf{sum}$  of $(E,\mathcal{F}_1),\dots,(E,\mathcal{F}_n)$
. In fact, the union of matroid need not be operated on the same ground set . Generally let $\mathcal{M}_i=(E_i,\mathcal{F}_i)$ , then we have the notion called \textbf{general union} in this paper;
\[\bigvee_{i\in J} \mathcal{M}_i:
=(\bigcup_{i\in J} E_i, \bigvee_{i\in J} \mathcal{F}_i)\],
where $\bigvee_{i\in J} \mathcal{F}_i$ is defined as 
\[\bigvee_{i\in J} \mathcal{F}_i=\{I\subseteq \bigcup_{i\in J} E_i:I=\bigcup_{i\in J} I_i,I_i\in\mathcal{F}_i\}\]
here $J$ is a finite index set .\\

\textit{Note.}  As standard notations we use $\mathcal{B}(\mathcal{M})$ to denote the set of all basis in a given matroid $\mathcal{M}$ and $\mathcal{I}(\mathcal{M})$ to denote the set of all independent sets in $\mathcal{M}$ .

\section{Proof of Main Results}\label{sec:DC}

Before proving the main results we are intend to propose a basic problem which is motivated by the set-theoretical  complementary relationship . In this paper , we say two independent systems $\mathcal{M}_1\subseteq \mathcal{M}_2$  which means that $\mathcal{I}(\mathcal{M}_1)\subseteq \mathcal{I}(\mathcal{M}_2)$ on a common ground set , where $\mathcal{I}(\mathcal{M})$ denotes the set consisting of independent sets in a given independent system $\mathcal{M}$ .

\begin{que}

\label{que}
Does the following equality 
\[(E,(\bigcap_{i=1}^n \mathcal{F}_i)^*)=(E,\bigsqcup_{i=1}^n\mathcal{F}_i^*)\]
hold or not ? If this equality does not hold generally , when does it hold ?
\end{que}

Generally this equality does not hold because $\bigl(E,(\bigcap_{i=1}^n \mathcal{F}_i)^*\bigr)$ is not a matroid generally . However , applying $(i)$ of \ref{dualrank} if the above equality holds , it implies that $(E,(\bigcap_{i=1}^n \mathcal{F}_i)^*)$ 
must be a matroid for the reason that  $(E,\bigsqcup_{i=1}^n\mathcal{F}_i^*)$ is a matroid which is easily to check to be  an basic property for the union of matroids .
Before giving the answer of this question , we need some set-theorectical preparations . 
\begin{lemma}
 Let $(E,\mathcal{F}_1),\dots,(E,\mathcal{F}_n)$ be as above .

If $X\subseteq E$ and $X=\bigcup_{i=1}^n X_i$ , where $X_i\in \mathcal{F}_i$ respectively , then X is partible , i.e. , 
$X=\bigsqcup_{i=1}^n {X_i}'$ , here ${X_i}' \subseteq X_i$ respectively for $i=1,\dots,n$ .
\end{lemma}

\begin{proof}
We will prove this lemma by induction . For the case $n=1$ the lemma is true automatically . Suppose that the property is true when $n=k$ , that is , for any set $X$ that is the union of k subsets of $X$ , explicitly $X=\bigcup_{i=1}^k X_i$ , then we have $X=\bigsqcup_{i=1}^k X_i'$ , where $X_i'\subseteq X_i$ respectively for $i=1,\dots,k$ . \\
Suppose that $X=\bigcup_{i=1}^{k+1} X_i$ , then 

\begin{equation}
    \begin{split}
        \bigcup_{i=1}^{k+1} X_i &=\biggl(\bigcup _{i=1}^{k} X_i\biggr) \bigcup X_{k+1}\\
        &=\biggl(\bigsqcup_{i=1}^{k} X_i'\biggr) \bigcup X_{k+1}\\
        &=\biggl(\bigsqcup_{i=1}^{k} X_i'\biggr) \bigsqcup    \Biggl( X_{k+1}\backslash \biggl(\bigsqcup_{i=1}^{k} X_i'\biggr)\Biggr)
    \end{split}
\end{equation}
where $X_i'\in X_i$ for $i=1,\dots,k$ by assumption . Moreover , because the subset of an independent set is also an independent set , we  complete the proof .

\end{proof}

\begin{lemma}
\label{dcon}
Let $(E,\mathcal{F}_1),\dots,(E,\mathcal{F}_n)$ be as above . Then we have the following containment relationship :
\[(\bigcap_{i=1}^n \mathcal{F}_i)^*\subseteq \bigsqcup_{i=1}^n\mathcal{F}_i^* .\]
\end{lemma}

\begin{proof}
For every $X\subseteq E$ ,  $X\in (\bigcap_{i=1}^n \mathcal{F}_i)^*\iff $ there is a basis $B_0\in\bigcap_{i=1}^n \mathcal{F}_i$  such that $X\bigcap B_0=\emptyset$ . Now we extend $B_0$ to $B_1,\dots,B_n$ respectively , which are  bases for $\mathcal{F}_1,\dots,\mathcal{F}_n$ respectively . On the one hand , 
\[(E\backslash B_0)=(E\backslash B_1)\bigcup\dots\bigcup(E\backslash B_n)\]
and it is obvious that $E\backslash B_i\in \mathcal{F}_i^*$ for $i=1,\dots,n$ . Nextly using $\mathbf{Lemma\;12}$ we obtain that $E\backslash B_0$ is partible , hence $E\backslash B_0\in\bigsqcup_{i=1}^n\mathcal{F}_i^*$ . On the other hand , for $(X\subseteq E)\backslash B_0$ , it results that $X\in\bigsqcup_{i=1}^n\mathcal{F}_i^*$ . This verifies the containment relation .
\end{proof}

Actually we have a stronger version of the following easy fact which enables us to deduce a nontrivial lower bound for the maximum cardinality of the common independent sets in the  intersection of several matroids .

\begin{proposition}
\label{prop3.1}
Let $\mathcal{M}_i(E,\mathcal{F}_i)$ be n matroids where $i=1,\dots,n$ . And set $\{T_i\}_{i=1}^n$ is a covering and co-covering of $E$ which we have explained in \textbf{Theorem}\ref{thm1.3}
. Futhermore we assume $T_i\in \mathcal{B}(\mathcal{M}_i)$ . Then we  have 
\[(\bigcap_{i=1}^n \mathcal{M}_i)^*\subseteq \bigvee_{i=1}^n (\mathcal{M}_i/T_i)^*\].
\end{proposition}

\begin{proof}
Let $\mathcal{B}_i$ denote the set of all bases in $\mathcal{M}_i$ and $\Bar{\mathcal{B}}$ the set of bases in $\bigcap_{i=1}^n \mathcal{M}_i$  . Then following the definition of the dual, we have 
\[(\bigcap_{i=1}^n \mathcal{M}_i)^*=\{C\subseteq E:C\cap B =\emptyset,B\in \Bar{\mathcal{B}}\}\]
.It suffices to show that every basis $C $ in $(\bigcap_{i=1}^n \mathcal{M}_i)^*$
can be covered by the union of n bases in $(\mathcal{M}_i/T_i)*$ respectively . There is a basis $B\in \Bar{\mathcal{B}}$ such that $C\cap B=\emptyset$ . On the one hand , by the assumptions on $\{T_i\}_{i=1}^n$ we can see that $C\cap (\cup_{i=1}^n T_i)=C$ and $C\cap (\cup_{i=1}^n T_i^c)=\cup_{i=1}^n(C\cap T_i^c)  =C$ . On the other hand , for  $T_i\in \mathcal{B}(\mathcal{M}_i)$ , it is evident that  $T_i^c\in \mathcal{B}(\mathcal{M}_i^*)$ . It follows that $C\cap T_i^c \in \mathcal{F}_i^* $ . Next, we can extend $C\cap T_i^c$ to a basis $B_i^*\in \mathcal{F}_i^*$ resp. for all $i$ . Therefore $C=\cup_{i=1}^n(B_i^*\setminus T_i)$.  This yields a covering of $C$ by n distinct bases in $\mathcal{M}_i^* \setminus T_i$ for $i=1,\dots,n$ .

And we know that $(\mathcal{M}_i/T_i)^*=\mathcal{M}_i^* \setminus T_i$ . This leads to the proof .

\end{proof}

\begin{proposition}
\label{dualtwo}
Let $(E,\mathcal{F}_1)$ and $(E,\mathcal{F}_2)$ be two matroids and assume their intersection is also a matroid . Then the dual relation in $\mathbf{Question}$ \ref{que} holds . That is :
  \[(E,(\bigcap_{i=1}^2 \mathcal{F}_i)^*)=(E,\bigsqcup_{i=1}^2\mathcal{F}_i^*) .\]
\end{proposition}

\begin{proof}
It is equivalent to prove that :
  \[(E,(\bigcap_{i=1}^2 \mathcal{F}_i))=(E,(\bigsqcup_{i=1}^2\mathcal{F}_i^*)^*)\]
  We just need to show that , $(1)$ for every $Y\in \bigcap_{i=1}^2 \mathcal{F}_i$ , there is a basis $C\in  \bigsqcup_{i=1}^2\mathcal{F}_i^*$ such that $Y\bigcap C=\emptyset$ ; $(2)$ for every $Y'\in(\bigsqcup_{i=1}^2\mathcal{F}_i^*)^*$ , $Y'$ is also in $\bigcap_{i=1}^n \mathcal{F}_i$ .\\

 To prove $(1)$ we just need to construct such $C$ .     $C$ canbe obtained in the following way  : since $Y\in\bigcap_{i=1}^2 \mathcal{F}_i$ we extend $Y$ to $B_0$  , which is a  basis for $\mathcal{F}_1\bigcap\mathcal{F}_2$ ( Since $\mathcal{F}_1\bigcap\mathcal{F}_2$ is a matroid , the following extension is well-defined ) .  We nextly find a minimum cardinality subset $E\supseteq B_m\supseteq B_0$ satisfying that $B_m=B_1\bigcup B_2$ where $B_1,B_2$ are bases for $\mathcal{F}_1,\mathcal{F}_2$ respectively and they both contain $B_m$ . Then we assert that $E\backslash B_m$ is a basis for $(\bigcap_{i=1}^2 \mathcal{F}_i)^*$ . \\
 
 Fisrtly $E\backslash B_m=E\backslash(B_1\bigcup B_2)=(E\backslash B_1)\bigcup(E\backslash B_2)$ is obviously in $\bigsqcup_{i=1}^2\mathcal{F}_i^*$ . Then we just need to prove that $E\backslash B_m$ is a maximal cardinality independent set in $\bigsqcup_{i=1}^2\mathcal{F}_i^*$ .If this is not true , there are another two bases $B_1 ',B_2 '$ of $\mathcal{F}_1,\mathcal{F}_2$ respectively such that $(B_1 '\bigcup B_2 ' )\bigcap (E\backslash B_m=E) =\emptyset $ and $|B_1 '\bigcup B_2 '|<|B_1 \bigcup B_2 |=|B_m|\iff |B_1 '\bigcap B_2 '|>|B_1 \bigcap B_2 |=|B_m|$ for the reason that $|B_1 '| =|B_1|$ and $|B_2 '| =|B_2|$ and the inclusion and exclusion principle . However $|B_1 \bigcap B_2 |=|B_m|$ is the maximal cardinality set in $\mathcal{F}_1\bigcap\mathcal{F}_2$ , this lead to a contradiction . Therefore $E\backslash B_m$ is a basis for $(\bigcap_{i=1}^2 \mathcal{F}_i)^*$ . Put $C=E\backslash B_m$ , we can finish the proof of $(1)$ .\\

  It remains to prove $(2)$ . By definition , there is a basis $D\in \bigsqcup_{i=1}^2\mathcal{F}_i^*$ , such that $Y'\bigcap D=\emptyset$ . We can write $D=D_1\bigsqcup D_2$ , here $D_i\in \mathcal{F}_i^*$ for $i=1,2$ respectively . Extend $D_1$ to a basis of $\mathcal{F}_1^*$  denoted as $\Bar{D}_1$ . However , since $D=D_1\bigsqcup D_2$ is a basis , adding any elements in $E\backslash D$ wouldn't change the rank of $D$ , that is to say , $\Bar{D}_1$ must contains in $D$
  no matter how we extend $D_1$ . On the other hand , since $Y'\bigcap D=\emptyset$ , it forces $Y'\bigcap \Bar{D}_1=\emptyset$ . We therefore can assert that $Y'\in \mathcal{F}_1$ . In the same way we can also prove that $Y'\in \mathcal{F}_2$  . Then we conclude that $Y'\in\bigcap_{i=1}^2 \mathcal{F}_i$ . This leads to the proof of $(2)$ .

\end{proof}

\begin{corollary}[Proposition\ref{prop1.2}]
Let $(E,\mathcal{F}_1),\dots,(E,\mathcal{F}_n)$ be the same as stated in \textbf{Proposition} \ref{prop1.1} .  And set $(E,\mathcal{F}_j^{*})$ the duality of $(E,\mathcal{F}_j)$ for $j=1,\dots,n$ . Moreover we assume that $(E,\bigcap_{i=1}^n\mathcal{F}_i)$ is an independent system with rank function $\Bar{r}$ ( Notice that the intersection of matroids are not generally a matroid ) . Then we have 
\begin{enumerate}
    \item \[\Bar{r}^*(X)\leqslant \min_{A \subseteq X}\biggl(|X|+(n-1)|A|+\sum_{i=1}^n r_i(E\backslash A)-\sum_{i=1}^n r_i(E)\biggr)\]
    Here $\Bar{r}^*$ is the rank function of $(E,(\bigcap_{i=1}^n\mathcal{F}_i)^*)$ . Especially if we take $X=E$ , we have the following :
  \[\max\{|X|:X\in \bigcap_{i=1}^n\mathcal{F}_i\}\geqslant \max_{A \subseteq E}\biggl(\sum_{i=1}^n r_i(E)-\sum_{i=1}^n r_i(E\backslash A)-(n-1)|A|\biggr).\]
  \item If we restrict ourselves on two matroids  $(E,\mathcal{F}_1),(E,\mathcal{F}_2)$ , and furthermore assume that $\mathcal{F}_1\bigcap \mathcal{F}_2$ is a matroid , then the two following inequalities are actually equalities which agree with $Edmonds'\mathbf{\;min-max \;Theorem}$ \ref{Edmonds} .
\end{enumerate}

\end{corollary}

\begin{proof}

Due to $\mathbf{Proposition\;}$ \ref{dcon} we have 
\[(\bigcap_{i=1}^n \mathcal{F}_i)^*\subseteq \bigsqcup_{i=1}^n\mathcal{F}_i^* .\]

therefore
\[\Bar{r}^*(X)\leqslant \hat{r}(X)\]
here $\hat{r}(X)$ is the rank function for $\bigsqcup_{i=1}^n\mathcal{F}_i^*$ . Applying $\mathbf{Theorem\;}$ \ref{Nash} and $\mathbf{Proposition\;}(ii)$ \ref{dualrank} we have 

\begin{equation}
    \begin{split}
        \Hat{r}(X) \\
     &= \min_{A\subseteq X} \bigl(|X\backslash A|+\sum_{i=1}^n r_i^*(A)\bigr) \\
              &= \min_{A\subseteq X}\bigl(|X\backslash A|+\sum_{i=1}^n(|A|+r_i(E\backslash A)-r_i(E)\bigr) \\
              &= \min_{A\subseteq X}\biggl((n-1)|A|+|X|+\sum_{i=1}^n\bigl(r_i(E\backslash A)-r_i(E)\bigr)\biggr)
    \end{split}
\end{equation}

hence 
\[\Bar{r}^*(X)\leqslant \min_{A \subseteq X}\biggl(|X|+(n-1)|A|+\sum_{i=1}^n r_i(E\backslash A)-\sum_{i=1}^n r_i(E)\biggr)\]

What's more , since
\\\\
\[\Bar{r}^*(X)=|X|+\Bar{r}(E\backslash X)-\Bar{r}(X)\]
Taking $X=E$ we obtain 
\\\\
\[|E|-\Bar{r}(E)=|E|-\max\{|X|:X\in \bigcap_{i=1}^n\mathcal{F}_i\}\leqslant \min_{A \subseteq E}\biggl(|E|+(n-1)|A|+\sum_{i=1}^n r_i(E\backslash A)-\sum_{i=1}^n r_i(E)\biggr)\]
\\\\
After eliminating $|E|$ in both side we get 
 \[\max\{|X|:X\in \bigcap_{i=1}^n\mathcal{F}_i\}\geqslant \max_{A \subseteq E}\biggl(\sum_{i=1}^n r_i(E)-\sum_{i=1}^n r_i(E\backslash A)-(n-1)|A|\biggr).\]
 
This proves part 1 . And part 2 is evident for the reason that both inequalities in part 1 can be strengthened to be equalities  according to $\mathbf{Proposition\;}$ \ref{dualtwo} . 
\end{proof}

\begin{corollary}[Theorem\ref{thm1.3}]
We omit the statement of \textbf{Theorem}\ref{thm1.3}  . For details the readers can check \textbf{Section 1} . 
\end{corollary}

\begin{proof}
The proof is similar with that of \textbf{Proposition}\ref{prop1.2} . Applying \textbf{Proposition}\ref{prop3.1} , we briefly write $\bigcap_{i=1}^n\mathcal{M}_i$ as $\mathcal{P}$ , write $\bigvee_{i=1}^n(\mathcal{M}_i/T_i)^*$ as $\mathcal{N}$ , then we have
\[r_{\mathcal{P}}^*(S)\leqslant r_{\mathcal{N}}(S)\]
for any $S\subseteq E$ . Notice that we cannot apply \textbf{Theorem}\ref{thm1.2} to  $r_{\mathcal{N}}(S)$ .However after extending the domain of  each matroid $\mathcal{M}_i/T_i$ to $E$ trivially , i.e. adding elemnts of $T_i$ to the ground set of  $\mathcal{M}_i/T_i$ acting as $|T_i|$ loops which means that these additional elements vanish in the rank functions . Then we can utilise \textbf{Theorem}\ref{thm1.2} to compute $r_{\mathcal{N}}(S)$ . And combining
\textbf{Proposition 2.4} \ref{prop2.4.3}) and \textbf{Proposition} \ref{dualrank} , \textbf{Theorem}\ref{thm1.3} canbe deduced after some straightforward computations .

\end{proof}

 The remaining part of this paper is left to prove \textbf{Proposition} \ref{prop1.1} and \textbf{corollary}\ref{cor1.1} .  Let's recall what $\mathbf{Theorem\;}$ \ref{prop1.1} tells us .
 
 \begin{theorem}[Theorem\;\ref{prop1.1}]
Let $(E,\mathcal{F}_1)$ , $(E,\mathcal{F}_2)$ , ... , $(E,\mathcal{F}_n)$ be n  matroids with rank functions $r_1,\dots,r_n$ respectively . And assume that $(E,\mathcal{F}_1\bigcap\mathcal{F}_2)$ , $(E,\mathcal{F}_1\bigcap\mathcal{F}_2\bigcap\mathcal{F}_3)$ , $\dots$ , $(E,\bigcap_{i=1}^n\mathcal{F}_i)$ are all matroids . And for any $X\subseteq E$  we take a increasing filtration of $X$ randomly :
\[X\subseteq X_1\subseteq X_2\subseteq \dots\subseteq X_{n-2}\subseteq E\]
Then we have the following min-max result :
\begin{equation}
    \begin{split}
        \max\{|S|:S\in\bigcap_{i=1}^n\mathcal{F}_i\}&=\min\{r_1(X)+r_2(X_1\backslash X)+\dots+r_{n-1}(X_{n-2}\backslash X_{n-3})+r_n(E\backslash X_{n-2}):\\ &\quad X\subseteq X_1\subseteq X_2\subseteq \dots\subseteq X_{n-2}\subseteq E\}
    \end{split}
\end{equation}
\end{theorem}

Before proving this theorem , we will introduce two new concepts : $\mathbf{submatroid}$ of a given matroid and $\mathbf{restriction}$ of a matroid to its subset . Set $M=(E,\mathcal{F})$ a matroid , here $E$ is a finite set .

\begin{definition}[Submatroid]

A submatroid $M_0$ of $M$ is a matroid $(S,\mathcal{F}_S)$ where $S\subseteq E$ is a subset of $E$ and $\mathcal{F}_S$ is defined as follow :
\[\mathcal{F}_S=\{F\in 2^{S}:F\in \mathcal{F}\}\]

\end{definition}

It is easy to check that $(S,\mathcal{F}_S)$ satisfies all the matroid axioms $(M_1),(M_2),(M_3)$ . Hence the notion of submatroid  is well defined .

\begin{definition}[Restriction of a matroid]
Let $M=(E,\mathcal{F})$ be a matroid as above with rank function $r$ . For every subset $C\subseteq E$ , the $\mathbf{restriction}$ of $M$ to $C$ is nothing but the submatroid $(C,\mathcal{F}_C)$ defined as above , for brevity we denote it as $M|_C$ or $(E,\mathcal{F})|_C$ . \\

\end{definition}

And futhermore we denote the rank function of $M|_C$ as $r|_C$ .

\begin{lemma}
\label{rest}
If $K\subseteq E$ and even satisfying that $K\subseteq C\subseteq E$ , then we have
\begin{equation}
\begin{split}
    r(K)&=r|_C(K).
\end{split}
\end{equation}
\end{lemma}
\begin{proof}
On the one hand ,  by definition of rank function obviously we have :
\[r(K)\geqslant r|_C(K)\]

On the other hand for every $K\subseteq C$ , let $P$ be the maximal cardinality independent set contained in $K$ . It is straightforward to see that $P\subseteq C$ , it implies that $P$ is also in $\mathcal{F}_C$ . Therefore 
\[r(K)\leqslant r|_C(K)\]

Finally we can conclude that
\[r(K)= r|_C(K).\]

\end{proof}

\begin{proof}[Proof\;of \;Proposition\;\ref{prop1.1}]
Our strategy  is to prove by induction . For $n=2$ , \textbf{Proposition}\ref{prop1.1} is nothing but a special case of  Edmonds' min-max Theorem $(\mathbf{Theorem\;}\ref{Edmonds})$ . \\

Assume that equality $(4)$ holds for $n\leqslant k$ , here $k$ is a positive integer bigger  than 2 . And set $\Bar{r}_k$ the rank function for $(E,\bigcap_{i=1}^k\mathcal{F}_i)$ . Due to $\mathbf{Theorem}$ \ref{Edmonds} we have 
\begin{equation}
    \begin{split}
        \max\{|S|:S\in \bigcap_{i=1}^{k+1}\mathcal{F}_i\} 
        &=\max\{|S|:S\in(\bigcap_{i=1}^{k}\mathcal{F}_i)\bigcap\mathcal{F}_{k+1}\} \\
        &=\min\{\Bar{r}_k(X)+r_{k+1}(E\backslash X):X\subseteq E\}
    \end{split}
\end{equation}
The second equality is a direct result of application of Edmonds'min-max Theorem for $\bigcap_{i=1}^{k}\mathcal{F}_i$ and $\bigcap_{i=1}^{k}\mathcal{F}_i$ .\\
And observe that $\Bar{r}_k(X)=\Bar{r}_k|_X(X)$ according to $\mathbf{Lemma\;}$ \ref{rest} . Here $\Bar{r}_k|_X$ is the rank function for a submatroid  $(X,(\bigcap_{i=1}^{k}\mathcal{F}_i)|_X)$ of $(E,\bigcap_{i=1}^{k}\mathcal{F}_i)$ . What's more , we notice that
\[(\bigcap_{i=1}^{k}\mathcal{F}_i)|_X=\bigcap_{i=1}^{k}(\mathcal{F}_i|_X)\]
Applying the inductive hypothesis  to $(X,\bigcap_{i=1}^{k}(\mathcal{F}_i|_X))$ , we have 
\begin{equation}
\begin{split}
      \Bar{r}_k(X)&=\Bar{r}_k|_X(X) \\
      &=\min\{r_1|_X(A)+r_2|_X(A_1\backslash A)+\dots+r_{k-1}|_X(A_{k-2}\backslash A_{k-3})+r_k|_X(X\backslash A_{k-2}):\\ &\quad A\subseteq A_1\subseteq A_2\subseteq \dots\subseteq A_{n-2}\subseteq X\} \\
    &=\min\{r_1(A)+r_2(A_1\backslash A)+\dots+r_{k-1}(A_{k-2}\backslash A_{k-3})+r_k(X\backslash A_{k-2}):\\ &\quad A\subseteq A_1\subseteq A_2\subseteq \dots\subseteq A_{n-2}\subseteq X\}
\end{split}
\end{equation}
Finally combining $(6)$ and $(7)$ leads to a proof .
\end{proof}

\begin{proof}[Proof\;of\;Corollary\ref{cor1.1} ]
We just need to show that for any partition $\bigsqcup_{i=1}^m Y_i=E$
and every $I\in\mathcal{B}(\bigcap_{j=1}^m\mathcal{M}_j)$ we have  

\begin{equation}
\label{eq14}
|I|\leqslant \sum_{i=1}^m r_i(Y_i)
\end{equation}
It's clear that $I=\bigsqcup_{i=1}^m I\cap Y_i$ and $r_i(Y_i)\geqslant |I\cap Y_i|$ , so the inequality\ref{eq14} follows from the above facts .

\end{proof}

\textbf{Acknowledgements: } The author thanks the University of Chinese Academy of Sciences and his classmates and friends for their selfless help when he wrote this paper .

\printbibliography[]

\end{document}